\documentclass[12pt]{amsart}
\usepackage{amsfonts,amsmath,graphicx}
\usepackage{amsfonts}
\usepackage{amssymb}
\usepackage{amsmath}   
\usepackage{amstext}
\usepackage{amsbsy}
\usepackage{amsopn}
 \usepackage{amsthm}
 \usepackage{amscd}
\usepackage{amsxtra}
\usepackage{euscript}
\parindent
0pt

\newtheorem{thm}{Theorem}[section]

\newtheorem{prop}[thm]{Proposition}

\newtheorem{defn}[thm]{Definition}
\newtheorem{exmp}[thm]{Example}
\newtheorem{rem}[thm]{Remark}

\def\P{\mathbb{P}}
\def\G{\mathbb{G}}

\def\O{{\mathcal O}}

\def\eqref#1{(\ref{#1})}
\def\:{\colon }
\begin{document}
\title{Geometry of syzygies via Poncelet varieties }
\author{Giovanna Ilardi, Paola Supino, Jean Vall\`es}

\address{Dip. Matematica e Applicazioni, Univ. di Napoli,
via Cinthia, 80126 Napoli (Italy);
 {email:  giovanna.ilardi@dma.unina.it } }
\address{Dip. Matematica, Univ. RomaTre,
Largo L. Murialdo, 00146 Roma (Italy);
 {email: psupino@uniroma3.it } }
\address{Laboratoire de Mathematiques appliques
de Pau et de Pays de l'Adour,
 Avenue de l'Universit\'e
 64000 Pau(France);
 {email:  jean.valles@univ-pau.fr } }
%
\maketitle
\begin{abstract}
 We consider the  Grassmannian $\mathbb{G}r(k,n)$ of $(k+1)$-dimensional linear
subspaces of
 $V_n=H^0({\P^1},\O_{\P^1}(n))$. We define  
  $\frak{X}_{k,r,d}$
as the classifying  space of the
$k$-dimensional linear systems  of degree $n$ on $\P^1$ whose
basis realize a fixed number of polynomial relations of fixed
degree, say a fixed number of syzygies of a certain degree.
The first result of this paper is the computation of the dimension of   $\frak{X}_{k,r,d}$.
In the second part we make a link between   $\frak{X}_{k,r,d}$ and the  Poncelet varieties.
 In particular, we prove that the existence of linear syzygies 
implies  the existence of singularities on the Poncelet varieties.
\end{abstract}%
%
\pagestyle{myheadings} \thispagestyle{plain} \markboth{\today}{}



\section{Introduction and set up}
In this paper we are interested in linear systems $\Lambda$ on $\P^1$ of degree 
 $n$ and  projective dimension $k$ (so, from now we assume that $n>k$), more particulary in those
having an  algebraic limitation, namely   the syzygies.
A  {\em syzygy} of degree $d$ for $\Lambda$ is a
$(k+1)$-uple of homogeneous forms of degree $d$, $(g_0,..,g_k)$,
 such that
$\sum_{i=0}^k g_if_i=0$. 
We say that $\Lambda$ has $r$ syzygies of degree $d$ if there
exist $r$ linearly independent $k+1$-uples $(g_{0,j},..,
g_{k,j})$, where $0\le  j \le r$.

  The locus
$\frak{X}_{k,r,d}$ lives inside $\mathbb{G}r(k,n)$ in a natural way:
\begin{align}\label{}
\frak{X}_{k,r,d}\colon &= \big\{ \Lambda\in \G r(k, n) \text { having } 
 \notag\\
 & \text {at least $r$ syzygies of degree $d$} \big\}. \notag
\end{align}

The first result of this paper is  the computation of the dimension of
$\frak{X}_{k,r,d}$.
The subvarieties
$\frak{X}_{k,r,d}$
turn out to be determinantal varieties for a suitable map of
vector bundles on the Grassmannian. This  extends  the main result (corollary 4.4) in \cite{IS}, where the computation was only proved
in the  $k=1$ case.
In a second part we  give a geometric interpretation of the varieties 
$\frak{X}_{k,r,d}$ in terms of Poncelet varieties. These varieties were 
introduced by Trautmann 
in \cite{T}, but, except for the case of curve, they have not been actually studied.
We prove that the existence of linear syzygies 
implies  the existence of singularities on the Poncelet varieties.

\section{The dimension of the varieties $\frak{X}_{k,r,d}$}
Let  $\Lambda$ be a linear system on $\P^1$ of degree 
 $n$ and dimension $k$. We choose
$u, v$ a system of coordinates on $\P^1$, and denote by $V_n$
the  $n+1$-dimensional vector space $H^0({\P^1},\O_{\P^1}(n))$ of
binary homogeneous forms of degree $n$, otherwise said binary {\em
quantics} of degree $n$. A base for $V_n$ is $x_0,\dots,x_n$ where
$x_i=u^iv^{n-i}$. \\
Choose a linear subspace $\Lambda$ of $\P V_n$, and let
$\{f_0,\dots,f_k\}$ be
 a basis for $\Lambda$.
 It   defines  a  morphism of vector bundles on
$\P^1$
\begin{equation}\label{phiLambda}
 \begin{CD}
 \phi_{\Lambda} :\, \Lambda\otimes \O_{\P^1} @>>> \O_{\P^1}(n)
\end{CD}
\end{equation}
which is surjective when $\Lambda$ has no base points.
The system $\Lambda$ gives a  map from $\P^1$ to $\P^k$
\begin{equation}\label{fLambda}
f_{\Lambda} :\, \P^1 \stackrel{\Lambda} \longrightarrow \P^k.
\end{equation}
Its image is a rational curve of degree
$n$ when $ \phi_{\Lambda}$ is surjective, and less than $n$ when $\Lambda$ has base points.
For general $\Lambda$, $\phi_{\Lambda}$  is a surjective morphism
of vector bundles on $\P^1$, thus there is  an exact sequence
\begin{equation}\label{ELambda}
 \begin{CD}
0 @>>> E_{\Lambda} @>>>\Lambda\otimes \O_{\P^1} @>>>
\O_{\P^1}(n)@>>> 0,
\end{CD}
\end{equation}
where the kernel is a vector bundle $E_{\Lambda}$ on  $\P^1$ of
rank $k$ and degree $-n$.\\

The short exact sequence
(\ref{ELambda}) twisted by $\O_{\P^1}(d)$
\[
 \begin{CD}
0@>>>E_{\Lambda}(d)@>>>\Lambda\otimes \O_{\P^1}(d) @>>>
\O_{\P^1}(n+d)@>>> 0
\end{CD}
\]
suggests that $\Lambda$ has exactly $r$ independent syzygies of
degree $d$ if and only if $h^0(E_{\Lambda}(d))=r$.
\\
Since any $k$-bundle of degree $-n$ splits as 
$\O_{\P^1}(-a_0)\oplus\dots\oplus
\O_{\P^1}(-a_k)$ for suitable positive $a_0\leq\dots\leq a_k$ such
that $a_0+\dots+a_k= n$, one can stratify the varieties $\frak{X}_{k,r,d}$ by all possible 
splitting  of the integer $n$ in $k$ pieces. This point of view is developed by Ramella  in \cite{R} to sudy  the stratification of the Hilbert scheme of rational curves $C$ embedded in projective space, by the splitting of the restriction of the tangent bundle to $C$, and by the splitting of the normal bundle. We will use this point of view 
(in thm \ref{dime}) in order to prove that the dimension is the expected one.
\begin{exmp}
When $n=5$ and $k=3$
the only possible cases for $E_{\Lambda}$ are:
\[
E_{\Lambda}=\O_{\P^1}(-1)\oplus \O_{\P^1}(-2)\oplus
\O_{\P^1}(-2),
\] which is the general case, and
\[
E_{\Lambda}=\O_{\P^1}(-1)\oplus \O_{\P^1}(-1)\oplus
\O_{\P^1}(-3).
\]
In the general case $h^0(E_{\Lambda}(1))=1$, in the other $h^0(E_{\Lambda}(1))=2$, thus the general stratum is 
$\frak{X}_{3,1,1}=\G r(3,5)$, and the stratum $\frak{X}_{3,2,1}$ is strictly contained in $\G r(3,5)$.
\end{exmp}
Note that if  $n<2k$ the general splitting has the first term $a_0=1$, thus $h^0(E_{\Lambda}(1))\geq 1$, 
which implies that  $\frak{X}_{k,1,1}^n= \G r(k,n)$, that
is: there always exists a linear syzygy.\\

\begin{thm}\label{dime}  $\mathrm{codim} (\frak{X}_{k,r,d},\G r(k,n))=(dk+k-n+r)r.$
Moreover the varieties $\frak{X}_{k,r,d}$ are Cohen-Macaulay with singular locus 
$\frak{X}_{k+1,r,d}$.
\end{thm}
\begin{proof}

Consider the
  universal vector bundle $\mathcal{U}=\mathcal{U}_{\G r(k,n)}$ on the Grassmannian
\[
\mathcal{U}=\{ (f,\Lambda)\in V_n\times \G r(k,n) \mid f\in\Lambda
\},
 \]
and the canonical map of vector bundles
$$
 \mathcal{U} \hookrightarrow V_n\otimes \O_{\G r(k,n)}
 $$
On  the product variety $\G r(k,n)\times \P^1$ let  $p$ and $q$ be
the two projections
\[
 \begin{CD}
\G r(k,n) @<p<<  \G r(k,n)\times \P^1 @>q>> \P^1.
\end{CD}
\]
 The  map  $p$ composed
with the evaluation map gives on $\G r(k,n)\times \P^1$
\[
 \begin{CD}
 p^{*}\mathcal{U}
  @>>> V_n\otimes \O_{\G r(k,n)\times \P^1}@>>> q^{*}
  \O_{\P^1}(n).
\end{CD}
\]
For all $d$ we also have a morphism
\[
 \begin{CD}
 p^{*}\mathcal{U}\otimes  q^{*} \O_{\P^1}(d)
  @>>> q^{*} \O_{\P^1}(n+d).
\end{CD}
\]
If we take now the direct image $p_*$ on the Grassmannian,we
obtain
\[
 \begin{CD}
 \mathcal{U} \otimes V_d
  @>>> V_{n+d}
\end{CD}
\]
which is just the relative version of our map $ \phi_{\Lambda}$ in (\ref{phiLambda})
twisted by $H^0(\O_{\P^1}(d))$
\[
\Lambda \otimes H^0(\O_{\P^1}(d)) \rightarrow H^0(\O_{\P^1}(n+d)).
\]
 Therefore
\[
\frak{X}_{k,r,d}=\lbrace \Lambda \in \G r(k,n) \mid
rk(\Phi_{\Lambda})\le (k+1)(d+1)-r\rbrace.
\]

 Applying  Thom-Porteous formula we compute 
the expected codimension for
$\frak{X}_{k,r,d}$ as   
$r(n+r-(d+1)k)$.
\\
It is a classical fact that when the codimension is exactly the expected one than the Chow class is
$det(c_{n-k(d+1)+r+j-i}(\mathbb{C}^{(d+1)(n+d+1)}\otimes\mathcal{U}^*))$.
\\
We  compute now the codimension of the tangent space in a generic
point ${\Lambda}$. Consider its associated bundle $E_{\Lambda}$, by  genericity we have
\[
E_{\Lambda}= \O_{\P^1}^r(-d) \oplus \O_{\P^1}^{k-r-B}(-A)
\oplus \O_{\P^1}^{B}(-A-1)
\]
 where $A$ and $B$ are uniquely
defined by
\[
(n-dr)=A(k-r)+B, \,\, 0\le B<k-r
\]
and, by
hypothesis on the syzygy, $A>d$. The
codimension of the tangent space in the point ${\Lambda}$ is then
$h^1(E_{\Lambda}\otimes E_{\Lambda}^{\vee})$.
 Since
\[
 E_{\Lambda}\otimes E_{\Lambda}^{\vee}=
\O_{\P^1}^{r(k-r-B)}(d-A) \oplus \O_{\P^1}^{rB}(d-A-1) \oplus R,
\]
where $R$ is a suitable bundle  with $h^1(R)=0$, we have $h^1(E_{\Lambda}\otimes E_{\Lambda}^{\vee})=r(n+r-(d+1)k)$   which is the expected codimension. 
\end{proof}

\section{Geometric description as Poncelet  varieties}
At the end 
of  his paper  \cite{T}, Trautmann has introduced a generalization of Poncelet curves, namely the Poncelet varieties, in higher dimension. Those are in bijective correspondence with the points of the Grassmannian.
The aim of this part is to describe the points of $\frak{X}_{k,r,d}$  as Poncelet varieties. In particular we will show that Poncelet varieties corresponding to  $\frak{X}_{k,1,1}$  are singular (see theorem \ref{teorema}).
Following \cite{V}, we define the Poncelet varieties  as determinant of sections of 
Schwarzenberger bundles, therefore we start by  recalling the definition of Schwarzenberger bundle and we 
describe the zero locus of their section ( see proposition \ref{prozero}). 
\subsection{Schwarzenberger bundles}
We denote by     $(x_i=u^{i}v^{k+1-i})$,
$(y_j=u^{j}v^{n-k-1-j})$, and $(z_{l}=u^{l}v^{n-l})$ the basis of
$V_{k+1}, V_{n-k-1}
\,\, \textrm{and} \,\,  V_{n}$ respectively. 
The multiplication of homogeneous polynomials in two variables
\[
 \begin{array}{c}
 V_{k+1}\otimes V_{n-k-1} \stackrel{\phi_{\times}} \longrightarrow V_{n},\\
(u^{i}v^{k+1-i}, u^{n-k-1-j}v^{j})\longmapsto   u^{n-j+i}v^{n-i+j}
\end{array}
\]
induces the following bundle on $\P V_n$:
\begin{equation}\label{SCH}
\begin{CD}
0 @>>> V_{n-k-1}\otimes \O_{\P V_{k+1}}(-1)@>A>>  V_{n}
\otimes \O_{\P V_{k+1}} @>B>>
E_{n}@>>> 0
\end{CD}
\end{equation}
where
$$A^T=
\left (
             \begin{array}{cccccccc}
              x_0 & x_1 & \cdots  & x_{k+1} & &   & &\\
               & x_0 & x_1 &\cdots      &x_{k+1} &  & &  \\
               &    &  x_0 & x_1 &\cdots      &x_{k+1} &  &   \\
               &   &  &\ddots & \ddots&  &\ddots&\\
                & & &   &x_0 & x_1 & \cdots  & x_{k+1}
             \end{array}
      \right )$$
and
$B=(z_0,\dots,z_{n})$. These bundles are called {\em Schwarzenberger bundles} because for $n=2$ they were first introduced by Schwarzenberger in \cite{S}.

Into the product $\P V_{k+1}\times\P V_{n}$ the projective bundle $\P E_{n}$ is defined by the equations 
$$
BA^T = ( \sum_{i=0}^{i=k+1}x_i z_{i+j}=0)_{j=0,\cdots,n-k-1} 
$$
\\
Let $s\in H^0(E_{n})=V_{n}$ be a non zero section.
We describe now the zero locus $Z(s)$ geometrically.
\\
We call $C_{k+1}^{\vee}\in \P (V_{k+1}^*)$ the rational normal curve of the $k$-osculating planes to $C_{k+1}$. 
\begin{prop}
\label{prozero}
Let   $s\in H^0(E_{n})$ be a non zero section and  $D_{n}(s)$
be the corresponding effective
divisor of degree $n$ on  $C_{k+1}^{\vee}$. We denote by  $Z(s)$ the
zero-scheme of $s$. Then,
$$a\in Z(s)   \Leftrightarrow a^{\vee}\cap C_{k+1}^{\vee} \subset
D_{n}(s).$$
 
More generally we have

$${\mathcal I}_{Z(s)}\subset {\frak m}_a^{r+1} \Leftrightarrow (a^{\vee})^{r+1}\cap
C_{k+1}^{\vee}\subset D_{n}(s).$$
\end{prop}
\begin{proof}
The section
$s$ corresponds to an
 hyperplane $H_s \subset \P(V_{n})$ or to an effective divisor of degree
$n$ on the
rational curve $C_{n}$, but also on $C_{k+1}^{\vee}$. The section $s$
$$
\begin{CD}
(E_{n})^{*} @>{s}>> {\O}_{\P V_{k+1}} @>>>
{\O}_{Z(s)}@>>>  0
\end{CD}
$$
induces a rational map
$
\P(V_{k+1}) \longrightarrow \P((E_{n})^{*})
$
which is not defined over the zero-scheme $Z(s)$.

We call $\pi$ the canonical projection $\pi:\P(E_{n})
\longrightarrow \P(V_{k+1})$.

Over a point
$a\in \P(V_{k+1})$ the fiber is 
$$\pi^{-1}(a)=\left\{(\sum_{i=0}^{i=k+1}a_iz_{i+j}=0)_{j=0,\dots,n-k-1}\right\}=
\P(E_{n}(a
))=\P^{k}.$$
The rational map $
\P(V_{k+1}) \longrightarrow \P((E_{n})^{*})
$
sends a point $a\in \P(V_{k+1}) $ onto $H_s\cap \pi^{-1}(a)$ which is, in
general, a $\P^{k-1}$, i.e a point in
$\P((E_{n})^{*}(a))=\P^{k\vee}$. This map is not defined when
$\pi^{-1}(a) \subset H_s $. The hyperplane
$H_s$ cuts the rational curve $C_{n}$ along an effective divisor
$D_{n}(s)$ of length $n$.  When $D_{n}(s)$ is smooth, there are $\binom{n}{k+1}$
subschemes $D_{k+1} \subset D_{n}(s)$ of length $k+1$; they generate
the $k$-planes $k+1$ secant to $C_{n}$ which are contained in
$H_s$.  Since $\P V_1 \simeq C_{n}\simeq C_{k+1}^{\vee}$ it is clear that  $D_{n}(s)$ corresponds to degree $n$ divisor on
$C_{k+1}^{\vee}$. We still denote it by $D_{n}(s)$. 

Then the zero-scheme $Z(s)$ is the set of points $a\in \P(V_{k+1})$
such that the divisor $a^{\vee}\cap
C_{k+1}^{\vee}$
 of degree $n$ belongs to $D_{n}(s)$. When $D_{n}(s)$ is smooth, we get $n$
osculating hyperplanes of
$C_{k+1}$ in $\P(V_{k+1})$. Every subset of $k+1$-osculating hyperplanes gives
a point in $\P(V_{k+1})$. These points
are   the zero-scheme of the section $s$. We have proved the first part of the 
proposition.
\\
Assume  ${\mathcal I}_{Z(s)}\subset {\frak
m}_a^{r+1}$. Let $H_s$ the hyperplane corresponding to the
section $s$ and $D_{k+1}$ the divisor on $C_{k+1}^{\vee}$ corresponding  to
$a$. Then $\mid(r+1)D_n\mid^{\vee}\subset H_s$. This proves
$(a^{\vee})^{r+1}\cap C_{k+1}^{\vee}\subset D_{n}(s)$.
\\
On the other hand, the inclusion  $\mid(r+1)D_{k+1}\mid^{\vee}\subset H_s$  proves that
the exceptional divisor
 $\mid D_{k+1} \mid^{\vee}$ of  $\P({\mathcal I}_{Z(s)})$ appears with  multiplicity
$(r+1)$, it means that  ${\mathcal I}_{Z(s)}\subset {\frak m}_a^{r+1}$.
\end{proof}
\subsection{Poncelet varieties}
The group  $SL(2, {\mathbb{C}})$   acts on $\G r(k,n)$ and we have an equivariant morphism $\label{Plu}\G r(k,n)\hookrightarrow \P (\bigwedge ^{k+1} V_n)$.
The $SL(2)-$modules $\bigwedge ^{k+1} V_n$ and  $S^{k+1}V_{n-k}$     are isomorphic (see \cite {FH}, p. 160).
Moreover by Hermite reciprocity formula (see \cite {FH}, p. 82 and p.160), we have $S^{k+1}V_{n-k}\cong S^{n-k}V_{k+1}$. So the Plucker embedding  becomes 

\begin{equation}\label{Tmap}
T\colon \G r(k,n)\hookrightarrow \P(S^{n-k}V_{k+1}).
\end{equation} 
It  associates to $\Lambda=\langle f_0,\dots,f_k\rangle$ the hypersurface of degree $n-k$ in $\P V_{k+1}$ with equation $f_0\wedge\dots\wedge f_k=0$.
Indeed, since from \eqref{SCH} 
\[
H^0(E_{n})\cong H^0(\O_{\P^1}(n)),\]
we consider $f_i$ as section of $E_n$.  We summarize these facts in
 the following commutative diagram:
\begin{equation}\label{comdia}
 \begin{CD}
@.@. \Lambda\otimes\O_{\P^{k+1}} @ =  \Lambda\otimes\O_{\P^{k+1}}\\
@. @. @VVV   @ VVV\\
 0@>>> V_{n-k-1}\otimes \O_{\P^{k+1}}(-1)@>>>V_n\otimes \O_{\P^{k+1}}  @>>> E_{n}@>>>0\\
  @.                           @V=VV              @ VVV         @VVV                         \\
0 @ >>> V_{n-k-1}\otimes \O_{\P^{k+1}}(-1)@>>>V_n/\Lambda\otimes \O_{\P^{k+1}}  @>>> {\mathcal L}@>>>0.
\end{CD}
\end{equation} 
where the support of the sheaf ${\mathcal L}$ is given by $f_0\wedge\dots\wedge f_k=0$.
\begin{defn} The 
  varieties defined as determinant of $k+1$ sections of  Schwarzenberger bundles 
on $\P^{k+1}$ are called {\em Poncelet varieties} of dimension $k$.
\end{defn} 
We point out that a section of $E_n$ corresponds to $n$ points on the rational normal curve $C_{k+1}$, and it vanishes along $\binom{n}{k+1}$
points in $\P^{k+1}$ which are the intersection points of $(k+1)$-osculating hyperplanes of $C_{k+1}$ in ${k+1}$ points chosen among the previous $n$ points. Hence these varieties are characterized by the nice following geometric fact: they contain the vertices of $k+1$ dimensional linear space of osculating polytopes. The case of curves is well-known since Darboux, but in higher dimension quite nothing exists in literature.

\begin{prop}
\label{Tpenc}
Let $\Lambda$  be a linear system in $\G r(k,n)$  with $d$ base points. Then $T(\Lambda)$
is the union of $d$ osculating hyperplanes to $C_{k+1}\subset \P V_{k+1}$ and a degree $(n-k-d)$ Poncelet variety of dimension $k$.
  \end{prop}
\begin{rem}In particular when $\Lambda$ is a pencil, a base point corresponds to a linear syzygy. The morphism $T$ sends  $
\frak{X}_{1,1,d}$ on the locus of Poncelet curves which are union
of a Poncelet curve of degree $d-1$ with
  $n-d$  tangent lines to  $C_2$.
\end{rem}
\begin{proof}
In the case of curves it is proved by Trautmann (see proposition 1.11 of \cite{T}).
In general assume that $\Lambda$ has $d<n$ base points and let $f=0$ be one equation for this base locus.
Then   we have a factorization:
\[ 
 \begin{CD}
\O_{\P^{1}}^{k+1}@>{\frac{\Lambda}{f}}>> \O_{\P^{1}}(n-d)\stackrel{f}\hookrightarrow  \O_{\P^{1}}(n).
\end{CD}
\]
The first arrow  gives the following vector bundle map:
\[ 
 \begin{CD}
 0@>>> \O^{k+1}_{\P(V_{k+1})}@>{\frac{\Lambda}{f}}>>E_{n-d}  @>>>K@>>>0
\end{CD},
\]
where $K$ is supported on a Poncelet variety of degree $n-k-d$.
The second arrow gives
\[ 
 \begin{CD}
 0@>>> E_{n-d}@>>>E_n  @>>> \oplus_{i=1}^{d}\O_{H_i}@>>>0
\end{CD}
\]
where $H_i$ is the osculating hyperplane to $C_{k+1}$ which corresponds to a base point on $C_{k+1}^{\vee}$.
Now the result follows from the following commutative diagram:
\[ 
 \begin{CD}
@. \O_{\P^{k+1}} @ =  \O_{\P^{k+1}}\\
@.  @V{\frac{\Lambda}{f}}VV   @ V{\Lambda}VV\\
 0@>>> E_{n-d}@>>>E_n  @>>> \oplus_ {i=1}^{d}\O_{H_i}@>>>0\\
  @.                           @VVV              @ VVV         @V=VV                         \\
0 @ >>> K@>>>{\mathcal L}  @>>> \oplus_{i=1}^{d} \O_{H_i}@>>>0,
\end{CD}
\]
where ${\mathcal L}$  is a sheaf supported by the Poncelet variety corresponding to $\Lambda$.
 \end{proof}

The following examples are done as suggested by the  commutative diagram  in  (\eqref{comdia}).

\begin{exmp}  In $\G r(2,4)$, $\frak{X}_{2,1,2}= \G r(2,4)$ and ${\mbox{\rm codim\,}} \frak{X}_{2,1,1}= 1$.

a) The net $\Lambda= \langle u^4,u^2v^2,v^4\rangle \in \frak{X}_{2,1,2}\setminus \frak{X}_{2,1,1}$ corresponds to the smooth Poncelet quadric 
$$ \mathrm{det} \left (
             \begin{array}{cc}
              x_1 & x_3   \\
              x_0&x_2
             \end{array}
      \right )=x_1x_2-x_0x_3 = 0.
$$

b) The net $\Lambda= \langle u^4,u^3v,v^4\rangle \in \frak{X}_{2,1,1}$  has associated Poncelet cone
$$ \mathrm{det} \left (
             \begin{array}{cc}
              x_2 & x_3   \\
              x_1&x_2
             \end{array}
      \right )=x_2^2-x_1x_3=0.
$$ 
\end{exmp}
\begin{exmp}\label{fixpt}
The net $ \Lambda= \langle u^3v,uv^3,v^4\rangle$ has a base point, its associated Poncelet quadric is
$$ \mathrm{det} \left (
             \begin{array}{cc}
              x_0 & x_2   \\
              0&x_1
             \end{array}
      \right )=x_1x_0=0.
$$
This quadric consists in  two planes, one of which osculating the rational normal curve $C_3$.
\end{exmp} 

\begin{thm}\label{teorema}
The Poncelet 
variety associated to any element of $\frak{X}_{k,1,1}$  is singular. Moreover the Poncelet variety associated to a general element of 
$\frak{X}_{k,1,1}$ contains $\binom{n-1}{k}$ lines and is singular in the  $\binom{n-1}{k+1}$ vertices of this configuration.
\end{thm}
\begin{proof}
Let $\Lambda$ have a syzygy of degree one : we can say that $\Lambda =\langle uf, vf,f_2,\dots,f_k\rangle$, where $f\in V_{n-1}$. 
The curve $\Gamma$ in $\P^{k+1}$ defined by the determinant of the two sections $uf, vf$ of $E_n$ is obtained as follows.
\\
The pencil $\langle uf,vf\rangle$  defines $n-1$ fixed points and a moving point $p$ on the rational normal curve $C_{k+1}^{\vee}$.
Therefore, the curve $\Gamma$ in $\P^{k+1}$ consists in $\binom{n-1}{k}$ lines. They are the
lines of $k$-planes in $\P^{(k+1)\vee}$ passing through $p$ and each subset of $k$ points chosen among the fixed ones.
\\
Each $(k+1)$-uple of points on $C_{k+1}^{\vee}$, chosen among the $(n-1)$ fixed one, gives a point on $\Gamma\subset \P^{(k+1)}$. It is the  intersection point
of $(k+1)$ lines corresponding to each choice of $k$ points among the $(k+1)$. Since the $(k+1)$ points are distinct on the rational normal curve $C_{k+1}^{\vee}$, the $\P^{k-1}$ generated by each choice of $k$ points do not have a common point, dually,  the configuration of lines is not contained in an hyperplane. 
Then, since the hypersurface defined by $\Lambda$ contains the curve $\Gamma$, it has singularities in the $\binom{n-1}{k+1}$ vertices of the configuration of lines $\Gamma$.
\end{proof}
  
\begin{exmp}
The net $\Lambda=\langle u^5+v^5,u^5-u^4v+u^3v^2-u^2v^3+uv^4,u^5-v^5\rangle$ has  a syzygy of degree $1$. It corresponds to 
the Poncelet cubic surface
$$ \mathrm{det} \left (
             \begin{array}{ccc}
              x_1+x_2 & x_2+x_3 &x_3  \\
              x_0+x_1&x_1+x_2&x_3+x_2\\
              x_0&x_1+x_0&x_1+x_2\\
             \end{array}
      \right )=0.
$$
This surface has four singular points, which is the maximum number of ordinary double points. Thus it is a Cayley cubic.
\end{exmp}
It could be interesting to explore the link between the sygyzies of higher degrees and the singularities 
of the associated Poncelet varieties.


\end{document}